\documentclass{article}

\usepackage{PRIMEarxiv}

\usepackage[utf8]{inputenc} % allow utf-8 input
\usepackage[T1]{fontenc}    % use 8-bit T1 fonts
\usepackage{hyperref}       % hyperlinks
\usepackage{url}            % simple URL typesetting
\usepackage{booktabs}       % professional-quality tables
\usepackage{amsfonts}       % blackboard math symbols
\usepackage{nicefrac}       % compact symbols for 1/2, etc.
\usepackage{microtype}      % microtypography
\usepackage{lipsum}
\usepackage{fancyhdr}       % header
\usepackage{graphicx}       % graphics
\graphicspath{{media/}}     % organize your images and other figures under media/ folder

\title{Stability of stochastic dynamic systems of a random structure with Markov switchings in the presence of concentration points
\thanks{\textit{\underline{Citation}}: 
\textbf{Lukashiv, T.; Malyk, I.V.; Chepeleva, M.; Nazarov,
P.V. Stability of stochastic dynamic systems of a random structure with Markov switchings in the presence of concentration points. Pages.... DOI:000000/11111.}} 
}

\author{
  Taras Lukashiv \\
  Multiomics Data Science Research Group, NORLUX Neuro-Oncology Laboratory \\
  Department of Cancer Research \\
  Luxembourg Institute of Health, L-1210, Luxembourg, Luxembourg\\
  Department of  Mathematical Problems of Control and Cybernetics,\\ 
  Yuriy Fedkovych Chernivtsi National University, 58000, Chernivtsi, Ukraine\\
  \texttt{t.lukashiv@gmail.com} \\
   \And
  Igor V. Malyk\\
  Department of  Mathematical Problems of Control and Cybernetics,\\ 
  Yuriy Fedkovych Chernivtsi National University, 58000, Chernivtsi, Ukraine\\
  \texttt{i.malyk@chnu.edu.ua} \\
  \And
  Maryna Chepeleva \\
  Multiomics Data Science Research Group,  Department of Cancer Research \\
  Luxembourg Institute of Health, L-1445, Strassen, Luxembourg\\
  Faculty of Science, Technology and Medicine, University of Luxembourg, \\
  2 av. de l'Université, L-4365, Esch-sur-Alzette, Luxembourg \\
  \texttt{maryna.chepeleva@lih.lu} \\
  \And
  Petr. V. Nazarov \\
  Multiomics Data Science Research Group,  Department of Cancer Research \\
  Luxembourg Institute of Health, L-1445, Strassen, Luxembourg\\
  \texttt{petr.nazarov@lih.lu} \\
}

\begin{document}
\maketitle

\begin{abstract}
This article aims to investigate sufficient conditions for the stability of stochastic differential equations with a random structure, particularly in contexts involving the presence of concentration points. The proof of asymptotic stability leverages the use of Lyapunov functions, supplemented by additional constraints on the magnitudes of jumps and jump times, as well as the Markov property of the system solutions. The findings are elucidated with an example, demonstrating both stable and unstable conditions of the system.
\end{abstract}

% keywords can be removed

\section{Introduction}
In the vast majority of works with jump changes of the trajectory, it is assumed that the distance between jumps is not less than some $\delta$, i.e. $\| t_k-t_{k-1} \|>\delta$. According to this assumption, only a finite number of jumps occur on a finite interval, which is an important condition for proving the stability and exponential boundedness of the solution. In this case, the conditions of existence, unity, and stability of the systems of stochastic differential equations with jumps are reduced to the corresponding statements for systems without jumps.

This work considers the case in which jumps can be concentrated at some point, which leads to the following relationship
$$
\lim\limits_{k\rightarrow\infty}t_k=t_\infty<\infty.
$$

In this case, the cumulative effect of jumps can lead to a lack of stability of the system. This effect can be illustrated by a simple example of an ordinary differential equation
$$
dx(t)=-x(t)dt
$$
with jumps
$$
x(t_k )=x(t_k-)(1+k^2 )
$$
at the points
$$
t_k=\frac{\alpha}{k}, \alpha>0.
$$
It can be easily concluded that 
$$
\lim\limits_{t\rightarrow\alpha-}\| x(t) \|=\infty
$$
for $x(0)\neq 0$. This simple example indicates that the size of the jump plays an important role in the presence of condensation points in the system. 

The models proposed in this work can be applied in the field of catastrophe theory. Random events with a frequency that exponentially increases over a finite time interval are considered there. At the same time, solving applied problems, one should consider the possibility of the collapse of the system in case of large disturbances and/or small intervals between disturbances $g$. The conditions of instability in such cases may have the following form 

$$
\| t_k-t_{k-1} \|<\delta_{min}
$$
for some $k$.

\section{Problem statement}\label{sec2}

On the probabilistic basis $(\Omega,\mathfrak{F},F,\mathbf{P})$ \cite{LM1}, \cite{LM2},  we consider a stochastic dynamic system of random structure given by a stochastic differential equation (SDE)
\begin{equation}\label{eq1}
dx(t)=a(t,\xi(t),x(t))dt+b(t,\xi(t),x(t))dw(t),~t\in \mathbb{R}_{+}\backslash K,
\end{equation}
with Markov switching
\begin{equation}\label{eq2}
\Delta x(t)=g(t_k -,\xi(t_k -),\eta_k,x(t_k -)),~~~t_k \in K=\{t_n \Uparrow\}, 
\end{equation}
and initial conditions

\begin{equation}\label{eq3}
x(0)=x_0\in \mathbb{R}^m,~\xi(0)=y\in\mathbf{Y}, \eta_{0}=h\in\mathbf{H}. 
\end{equation}

Here $\xi(t), t\ge 0,$ is a Markov chain with a finite number of states $\mathbf{Y}=\{1,2,...,N_{\xi} \}$ and generator $Q=\{\tilde{q}_{ij}\}, i,j=\{1,...,N_{\xi} \}$; $\{\eta_k, k\geq0\}$ is a Markov chain with values in space $\mathbf{H}$ and with a transition probability matrix $\mathbb{P}_H$; $x:\left[0,+\infty\right)\times\Omega \rightarrow\mathbb{R}^m$; $w(t), t\ge0,$ is a $m$-dimensional standard Wiener process; the processes $w, \xi$ and $\eta$ are independent random processes \cite{LM1}, \cite{LM2}.

We denote by
$$
\mathfrak{F}_{t_k}=\sigma(\xi(s),w(s),\eta_e, s\le t_k, t_e\le t_k)
$$
the minimal $\sigma$-algebra with respect to which  $\xi(t), \forall t\in [0, t_k]$ and $\eta_n, n\le k$ are measured.

Measured by a set of variables functions $a:\mathbb{R}_{+}\times\mathbf{Y}\times\mathbb{R}^m\rightarrow\mathbb{R}^m$, $b:\mathbb{R}_{+}\times\mathbf{Y}\times\mathbb{R}^m\rightarrow\mathbb{R}^m\times\mathbb{R}^m$,  $g:\mathbb{R}_{+}\times\mathbf{Y}\times\mathbf{H}\times\mathbb{R}^m\rightarrow\mathbb{R}^m$ satisfy the boundedness condition and the Lipschitz condition

\begin{equation}\label{eq4}
\| a(t,y,x)\|^2+\| b(t,y,x)\|^2+\| g(t,y,h,x)\|^2\le C(1+\| x\|^2);
\end{equation}

\begin{equation}\label{eq5}
\| a(t,y,x_1)-a(t,y,x_2)\|^2+ \| b(t,y,x_1)-b(t,y,x_2)\|^2\le L\| x_1-x_2\|^2, x_1, x_2 \in \mathbb{R}^m;
\end{equation}

\begin{equation}\label{eq6}
\| g(t_k,y,h,x_1)-g(t_k,y,h,x_2)\|^2\le L_k\| x_1-x_2\|^2, x_1, x_2 \in \mathbb{R}^m, \sum\limits_{k=1}^{\infty}L_k<\infty. 
\end{equation}

Consider the case of a point of concentration of jumps, i.e. 
\[\lim\limits_{n\rightarrow\infty}t_n=t^{*}\in [0,T].\] 
Let's assume that the following relations are true:

\begin{equation}\label{eq7}
    \sum\limits_{k=1}^{\infty}\gamma_k<\infty, \gamma_k=\sup\limits_{x\in \mathbb{R}^m,~y\in\mathbf{Y}, h\in\mathbf{H}} \| g(t_k,y,h,x)\|
\end{equation}
and

\begin{equation}\label{eq8}
\lim\limits_{\varepsilon\downarrow0}\left(\ln\varepsilon+N_\varepsilon\sum\limits_{k=1}^{N_\varepsilon} L_k\right)=-\infty, N_\varepsilon:=inf\left\{k\ge1:\sum\limits_{m=k}^{\infty}\gamma_m<\varepsilon\right\}.
\end{equation}

The conditions (4)-(8) guarantee the existence of a strong solution to the Cauchy problem (1)-(3) \cite{LM3}.

We denote by 
$$
{\it \mathbf{P}}_{k} ((y, h, x), \Gamma \times G\times \mathbf{C}):=
$$
$$
:=\mathbf{P}((\xi(t_{k+1}),\eta_{k+1},x(t_{k+1}))\in\Gamma \times G\times \mathbf{C}\vert(\xi(t_{k}),\eta_{k},x(t_{k}))=(y,h,x))
$$
the transition probability of the Markov chain $(\xi (t_{k} ), \eta _{k}, x(t_k) )$, that determine the solution to the problem (1)-(3)  on the $k$-th step. 

\textbf{Definition 1.} Discrete Lyapunov operator $(lv_{k} )(y,h,x)$ on a sequence of measurable scalar functions $v_{k} (y, h, x){\it :\mathbf{Y}}\times {\it \mathbf{H}}\times {\it \mathbb{R}}^{m} \to {\it \mathbb{R}}^{1} , k\in {\it \mathbb{N}}\cup \{ 0\}, $ for the SDE (1) with Markov switchings (2) is defined by the equality
\[(lv_{k} )(y, h, x):= \]
\begin{equation} \label{eq9}
:= \int\limits_{{\it \mathbf{Y}}\times {\it \mathbf{H}}\times {\it \mathbb{R}}^{m} }{\it \mathbf{P}}_{k} ((y, h, x) (du\times dz\times dl))v_{k+1} (u, z, l)-v_{k} (y, h, x).
\end{equation}

Here $v_{k} (y, h, x), k\in {\it \mathbb{N}}$, is a Lyapunov function defined by the following definition.

\textbf{Definition 2.} The Lyapunov function for the system (1)-(3) is a sequence of non-negative functions $\left\{v_{k} (y, h, x),k\ge 0\right\},$ for whom

\begin{enumerate}
\item for all $k\ge 0, y\in {\it \mathbf{Y}}, h\in {\it \mathbf{H}}, x\in {\it \mathbb{R}}^{m} $ the discrete Lyapunov operator $(lv_{k} )(y, h, x)$ (9) is defined;

\item if $r\to \infty $
$$
 \bar{v}(r)\equiv \mathop{\inf  }\limits_{k\in {\it \mathbb{N}, y\in {\it \mathbf{Y}},}  {h\in {\it \mathbf{H}}, \| x\|\ge r} } v_{k}(y, h, x)\to +\infty ;
$$

\item if $r\to 0$
$$
 \underline{v}(r)\equiv \mathop{\sup  }\limits_{ {k\in {\it \mathbb{N}}, y\in {\it \mathbf{Y}},}  {h\in {\it \mathbf{H}}, \| x\|\le r} } v_{k}(y, h, x)\to 0;
$$
\end{enumerate}
where $\bar{v}(r)$ and $\underline{v}(r)$ are continuous and monotonic for $r>0$.

\textbf{Definition 3.} A system with a random structure (1)-(3) is called:

-- stable in probability, if for $\forall \varepsilon _{1} >0, \varepsilon _{2} >0$ it can specify $\delta >0$ such that the inequality $\| x\|<\delta $ implies the inequality
\begin{equation} \label{eq10}
{\it \mathbf{P}}\left\{\mathop{\sup }\limits_{t\ge 0 } \| x(t)\|>\varepsilon _{1} \right\}<\varepsilon _{2}
\end{equation}
for all $y\in {\it \mathbf{Y}}, h\in {\it \mathbf{H}}$;

-- asymptotically stochastically stable, if it is stable in probability and for any $\varepsilon >0$ exists $\delta _{2} >0$ such that
\begin{equation} \label{eq11}
\mathop{\lim }\limits_{T\to \infty } {\it \mathbf{P}}\left\{\mathop{\sup }\limits_{t\ge T} \| x(t)\|>\varepsilon \right\}=0
\end{equation}
for all $\| x\|<\delta _{2} $ , $y\in {\it \mathbf{Y}}, h\in {\it \mathbf{H}}$ and $T\ge 0$.

\textbf{Definition 4.} A system with a random structure (1)-(3) is called:

--  mean square stable, if for $\forall \varepsilon >0$ it can specify the following $\delta >0$, that the inequality $\| x\|<\delta$ implies the inequality
\begin{equation} \label{eq12}
{\it \mathbf{E}}\| x(t)\|^{2} <\varepsilon
\end{equation}
for all $t \in  [0,T], y\in {\it \mathbf{Y}}, h\in {\it \mathbf{H}}$;

-- mean square asymptotically stable, if it is mean square stable for any $T>0$ and
\begin{equation} \label{eq13}
\mathop{\lim }\limits_{t\to \infty } \mathop{\sup }\limits_{y\in \mathbf{Y}, h\in \mathbf{H}} {\it \mathbf{E}}\| x(t)\|^{2} =0.
\end{equation}

If (10)-(13) hold true for all $x\in \mathbb{R}^m$, then the system is stable in the corresponding probabilistic sense on the whole.

For solving the problem (1)-(3) on the intervals $[t_k,t_{k+1})$, the following estimate is obtained.

{\textbf{Theorem 1.}}
Let the coefficients $a, b$ of the equation (1) satisfy the condition of uniform boundedness (4), and the condition (6) holds for the function $g$.

Then for all $k\ge 0$ for a strong solution of the Cauchy problem (1)-(3) holds the next inequality 

\begin{equation} \label{eq14}
{\it \mathbf{E}}\left\{\mathop{\sup }\limits_{t_{k} \le t\leq t_{k+1} } \| x(t)\| ^{2} \right\}\le 9e^{5C}(1+2L_{k+1})\left[\mathbb{E}\| x(t_{k} )\| ^{2} +5C (t_{k+1} -t_{k} )\right].
\end{equation}

\textit{Proof of Theorem 1.}

We use the same methodology as in \cite{LM4}, \cite{LM5}. A strong solution of the Cauchy problem (1), (3) for all $t\in[t_k,t_{k+1}), k\geq0$, can be written in the integral form

\begin{equation}\label{eq15}
x(t)=x(t_{k} )+\int\limits_{t_{k} }^{t}a(\tau ,\xi (\tau ),x(\tau))d\tau+\int\limits_{t_{k} }^{t}b(\tau ,\xi (\tau ),x(\tau))dw(\tau).
\end{equation}

After squaring the left and right sides of (15), calculating $\sup$, and applying the Cauchy–Schwarz inequality, we obtain:
$$
\sup\limits_{t_k\leq t<t_{k+1}}\| x(t)\|^2\leq
$$
$$
\leq 3\cdot\sup\limits_{t_k\leq t<t_{k+1}}\left\{  \| x(t_k)\|^2+\|\int _{t_{k} }^{t}a(\tau ,\xi (\tau ),x(\tau ))d\tau\|^2+\right.
$$
$$
+\left.\| \int _{t_{k} }^{t}b(\tau ,\xi (\tau ),x(\tau ))dw(\tau )\|^2\right\}\le
$$

$$
\leq3\left[\sup\limits_{t_k\leq t<t_{k+1}}\| x(t_k)\|^2+\sup\limits_{t_k\leq t<t_{k+1}}\int\limits _{t_{k} }^{t}\| a(\tau,\xi(\tau),x(\tau))\|^2 d\tau +\right.
$$
$$
+\left.\sup\limits_{t_k\leq t<t_{k+1}} \| \int\limits _{t_{k} }^{t}\| a(\tau,\xi(\tau),x(\tau)) dw(\tau)\|^2\right].
$$

To the last inequality, we apply the conditional mathematical expectation operation with respect to the $\sigma$-algebra $\mathfrak{F}_{t_k}$ and, taking into account the properties of the Ito integral, we obtain

$$
\mathbb{E}\left\{\sup\limits_{t_k\leq t<t_{k+1}}\| x(t)\|^2/\mathfrak{F}_{t_k}\right\}\leq 3\left[\mathbb{E}\| x(t_k)\|^2+C(t_{k+1}-t_k)+C\int\limits _{t_{k} }^{t_{k+1}}\mathbb{E}\| x(\tau)\|^2d\tau+\right.
$$
$$
+\left. 4C\int\limits _{t_{k} }^{t_{k+1}}\mathbb{E}\| x(\tau)\|^2d\tau\right]=
$$
$$
=3\left[\mathbb{E}\| x(t_k)\|^2+5C(t_{k+1}-t_k)+5C\int\limits _{t_{k} }^{t_{k+1}}\mathbb{E}\| x(\tau)\|^2d\tau      \right].
$$
Using the Gronwall inequality, we obtain an estimate of
$$
\mathbb{E}\left\{\sup\limits_{t_k\leq t<t_{k+1}}\| x(t)\|^2/\mathfrak{F}_{t_k}\right\}\leq
$$
$$
\leq 3\left[\mathbb{E}\| x(t_{k} )\|^{2} +5C(t_{k+1} -t_{k} )\right]e^{5C} .
$$

For $t=t_{k+1}$ the strong solution of the system (1)-(3), obviously, must satisfy the inequality
$$
\mathbb{E}\left\{ \| x(t_{k+1} )\|^{2} /\mathfrak{F}_{t_k}\right\}\le 3 \left[ \mathbb{E}\left\{ \| x(t_{k+1} -)\|^{2} /\mathfrak{F}_{t_k}\right\} + \right.
$$
$$
+2\mathbb{E}\left\{ \| g(t_{k+1}-,\xi(t_{k+1}-),\eta_{k+1},x(t_{k+1}-))-g(t_{k+1}-,\xi(t_{k+1}-),\eta_{k+1},0) \|^{2} /\mathfrak{F}_{t_k}\right\} +
$$
$$
\left.+2\mathbb{E}\left\{ \| g(t_{k+1}-,\xi(t_{k+1}-),\eta_{k+1},0) \|^{2} /\mathfrak{F}_{t_k}\right\} \right]\leq
$$
$$
\leq 3\left[(1+2L_{k+1})\mathbb{E}\left\{\sup\limits_{t_k\leq t\leq t_{k+1}}\| x(t)\|^2/\mathfrak{F}_{t_k}\right\}+C\right].
$$

Combining the last two inequalities, we get the desired estimate (14).

\textit{End of proof of Theorem 1.}

\textbf{Remark 1.} We will consider the stability of the trivial solution $x\equiv 0$, i.e. the satisfying of (4), if $C=0$ \cite{LM6}, \cite{LM7}, \cite{LM8}.

\textbf{Remark 2.} Note that the Lipschitz condition (5) was not used in the proof of the theorem, i.e. any (not necessarily unique) solution to the problem (1)-(3) satisfies the condition of the theorem.

{\textbf{Theorem 2.}} Let:

1) the conditions (4)-(6) are met;

2)  the Lyapunov functions $v_{k} (y, h, x)$ and $a_{k} (y, h, x), k\ge 0,$ exist, such that, based on the system, the following inequality
\begin{equation} \label{eq16}
(lv_{k} )(y, h, x(t))\le -a_{k} (y, h, x(t)), k\ge 0,
\end{equation}
is correct.

Then the system of random structure (1)-(3) is asymptotically stochastically stable on the whole.

\textit{Proof of Theorem 2.}
 Define by ${ \mathfrak{F}}_{t_{k} }=\sigma (\xi(s),\eta_{e}, s\leq t_k, t_e\leq t_k) $ a minimal $\sigma $-algebra, relative to which are measured $\xi (t)$ for all $t\in [0 , t_{k} ]$ and $\eta _{n} $ for $n\le k$. The conditional mathematical expectation is calculated by the formula
$$
{\it \mathbf{E}}\left\{{v_{k+1} (\xi (t_{k+1} ), \eta _{k+1} , x(t_{k+1} )) \mathord{\left/{\vphantom{v_{k+1} (\xi (t_{k+1} ), \eta _{k+1} , x(t_{k+1} )) {\mathfrak{F}}_{t_{k} } }}\right.\kern-\nulldelimiterspace} {\mathfrak{F}}_{t_{k} } } \right\}=
$$
\begin{equation} \label{eq17}
=\int _{{\it \mathbf{Y}}\times {\it \mathbf{H}}\times {\it \mathbb{R}}^{m} }{\it \mathbf{P}}_{k} ((\xi (t_{k} ),\eta _{k}, x)(du\times dz\times dl)v_{k+1} (u, z, l)) .
\end{equation}

Then, by the definition of the discrete Lyapunov operator $(lv_{k} )(y, h, x)$ (see (9)) from equality (17), considering (16), we get the inequality
$$
{\it \mathbf{E}}\left\{{v_{k+1} (\xi (t_{k+1} ), \eta _{k+1} , x(t_{k+1} )) \mathord{\left/{\vphantom{v_{k+1} (\xi (t_{k+1} ), \eta _{k+1} , x(t_{k+1} )) \mathfrak{F}_{t_{k} } }}\right.\kern-\nulldelimiterspace} \mathfrak{F}_{t_{k} } } \right\}=v_{k} (\xi (t_{k} ), \eta _{k} , x(t_{k} ))+
$$
\begin{equation} \label{eq18}
 +(lv_{k} )(\xi (t_{k} ), \eta _{k} , x(t_{k} ))\le \bar{v}(\| x(t_{k} )\| ).
\end{equation}

From Theorem 1 (because the existence of the second moment implies the existence of the first moment) and from properties of the function $\bar{v}$ follows the existence of a conditional mathematical expectation of the left-hand side of the inequality (18).

Now, using(17), (18), we write the discrete Lyapunov's operator $(lv_{k} )(y, h, x)$, which given on the solutions (1)-(3):
$$
lv_{k} (\xi (t_{k} ), \eta _{k} , x(t_{k} ))={\it \mathbf{E}}\left\{{v_{k+1} (\xi (t_{k+1} ), \eta _{k+1} , x(t_{k+1} )) \mathord{\left/{\vphantom{v_{k+1} (\xi (t_{k+1} ), \eta _{k+1} , x(t_{k+1} )) \mathfrak{F}_{t_{k} } }}\right.\kern-\nulldelimiterspace} \mathfrak{F}_{t_{k} } } \right\}-
$$
\begin{equation} \label{eq19}
-v_{k} (\xi (t_{k} ), \eta _{k} , x(t_{k} ))\le -a_{k} (\xi (t_{k} ),\eta_k,x(t_k))\le 0.
\end{equation}

Then, at $k\ge 0$ the next inequality holds
\[{\it \mathbf{E}}\left\{{v_{k+1} (\xi (t_{k+1} ), \eta _{k+1} , x(t_{k+1} )) \mathord{\left/{\vphantom{v_{k+1} (\xi (t_{k+1} ), \eta _{k+1} , x(t_{k+1} )) \mathfrak{F}_{t_{k} } }}\right.\kern-\nulldelimiterspace} \mathfrak{F}_{t_{k} } } \right\}\le v_{k} (\xi (t_{k} ), \eta _{k} , x(t_{k} )).\]

This means that a sequence of random variables
$$v_{k} (\xi (t_{k} ), \eta _{k} , x(t_{k} ))$$
forms a supermartingale in relation to $\mathfrak{F}_{t_{k} } $ \cite{LM9}.

Taking the mathematical expectation of both parts of inequality (19), we summarize the obtained expressions for $k$ from $n\ge 0$ to $N$, and obviously, we have the next inequality:
$$
{\it \mathbf{E}}\left\{v_{N+1} (\xi (t_{N+1} ), \eta _{N+1} , x(t_{N+1} ))\right\} -{\it \mathbf{E}}\left\{v_{n} (\xi (t_{n} ), \eta _{n} , x(t_{n} ))\right\} =
$$
$$
=\sum _{k=n}^{N}{\it \mathbf{E}}\left\{lv_{k} (\xi (t_{k} ), \eta _{k} , x(t_{k} ))\right\}\le 
$$
\begin{equation} \label{eq20}
\le -\sum _{k=n}^{N}{\it \mathbf{E}}\left\{a_{k} (\xi (t_{k} ), \eta _{k} , x(t_{k} ))\right\}\le 0.
\end{equation}

Since a random variable $\mathop{\sup }\limits_{t_{k} \le t\leq t_{k+1} } \| x(t)\| ^{2} $ does not depend on events of $\sigma $-algebra $\mathfrak{F}_{t_{k} } $ \cite{LM10}, then
\begin{equation} \label{eq21}
{\it \mathbf{E}}\left\{{\mathop{\sup }\limits_{t_{k} \le t\leq t_{k+1} } \| x(t)\| ^{2} \mathord{\left/{\vphantom{\mathop{\sup }\limits_{t_{k} \le t\leq t_{k+1} } \| x(t)\| ^{2} \mathfrak{F}_{t_{k} } }}\right.\kern-\nulldelimiterspace} \mathfrak{F}_{t_{k} } } \right\}={ E}\left\{\mathop{\sup }\limits_{t_{k} \le t\leq t_{k+1} } \| x(t)\| ^{2} \right\},
\end{equation}
that is, the inequality (14) also holds for the simple mathematical expectation
$$
{\it \mathbf{E}}\left\{\mathop{\sup }\limits_{t_{k} \le t\leq t_{k+1} } \| x(t)\| ^{2} \right\}\leq 3\mathbb{E}\| x\| ^{2}. 
$$

Next, we have
\[ {\it \mathbf{P}}\left\{\mathop{\sup }\limits_{t\ge 0} \| x(t)\| >\varepsilon _{1} \right\}=\]
 \[ ={\it \mathbf{P}}\left\{\mathop{\sup }\limits_{n\in \mathbb{N}} \mathop{\sup }\limits_{t_{n-1} \le t\leq  t_{n}}  \| x(t)\| >\varepsilon _{1} \right\}\le \]
 \[ \le {\it \mathbf{P}}\left\{\mathop{\sup }\limits_{n\in \mathbb{N}} 3 \| x( t_{n-1})\| >\varepsilon _{1} \right\}\le \]
  \[ \le {\it \mathbf{P}}\left\{\mathop{\sup }\limits_{n\in \mathbb{N}}  \| x( t_{n-1})\| >\frac{\varepsilon _{1}}{3} \right\}\le \]
\begin{equation} \label{eq22}
\le {\it \mathbf{P}}\left\{\mathop{\sup }\limits_{n\in \mathbb{N}} v_{n-1} (\xi (t_{n-1} ),\eta _{n-1} , x(t_{n-1} )) \ge \bar{v}(\frac{\varepsilon _{1}}{3}  )\right\}.
\end{equation}

If $\sup \| x(t_{k} )\| \ge r$, then, based on the definition of the Lyapunov function, the next inequality holds:
\begin{equation} \label{eq23}
 \mathop{\sup }\limits_{k\ge 0 } v_{k} (\xi (t_{k} ),\eta _{k} , x(t_{k} ))\ge \inf\limits_{k\ge 0 , y\in {\it \mathbf{Y}}, {h\in {\it \mathbf{H}}, \| x\| \ge r} } v_{k} (y, h,x) =\bar{v}(r) .
\end{equation}

Now let's use the well-known inequality for nonnegative supermartingales \cite{LM1}, \cite{LM9} to evaluate the right-hand side of (22):
\[{\it \mathbf{P}}\left\{\mathop{\sup }\limits_{n\in \mathbb{N}} v_{n-1} (\xi (t_{n-1} ),\eta _{n-1} , x(t_{n-1} ))\ge \bar{v}(\frac{\varepsilon _{1}}{3} )\right\}\le \]
\begin{equation} \label{eq24}
\le \frac{1}{\bar{v}(\frac{\varepsilon _{1}}{3}  )} v_{k } (y, h, x)\le \frac{\bar{v}(\| x\| )}{\bar{v}(\frac{\varepsilon _{1}}{3} )} .
\end{equation}

Given inequality (22), inequality (24) make it possible to guarantee the fulfillment of inequality (10) of stability in probability on the whole of the system (1)-(3).

From the inequality (20) follows the estimate
$$
{\it \mathbf{E}}\{ v_{N+1} (\xi (t_{N+1} ),\eta _{N+1} , x(t_{N+1} ))\} \le v_{0} (y, h, x)-
$$
\begin{equation} \label{eq25}
 -\sum _{k=0 }^{N}{\it \mathbf{E}}\{ a_{k} (\xi (t_{k} ),\eta _{k} , x(t_{k} ))\} \le v_{0 } (y, h, x),
\end{equation}
for all $N\ge 0, y\in {\it \mathbf{Y}}, h\in {\it \mathbf{H}}, x\in {\it \mathbb{R}}^{m} $.

Since the sequence $\{ a_{k} \}, k\ge 0$ forms Lyapunov functions, there must exist continuous strictly monotone functions $\underline{a}(r)$ and $\bar{a}(r)$, which are zero if $r=0$ \cite{LM11} and such that
\begin{equation} \label{eq26}
\bar{a}(\| x\| )\le a_{k} (y, h, x)\le \underline{a}(\| x\| )
\end{equation}
for $\forall k\in \mathbb{N}, y\in {\it \mathbf{Y}}, h\in {\it \mathbf{H}}$ and $x\in {\it \mathbb{R}}^{m} $.

Thus, from the convergence of the series on the left side of the inequality (25) (which will be convergent in the case of convergence of the series $\sum\limits_{k=1}^{\infty}L_k$) follows the convergence of the series  $\sum\limits_{k=0 }^{\infty }{\it \mathbf{E}} \{ \bar{a}(\| x(t)\| \} $ for $\forall t \ge t_k, y\in {\it \mathbf{Y}}, h\in {\it \mathbf{H}}, x\in {\it \mathbb{R}}^{m} $.

Then, taking into account the continuity of $\underline{a}(r)$ and the equality $\underline{a}(0)=0$, we have:
\begin{equation} \label{eq27}
\mathop{\lim }\limits_{k\to \infty } \| x(t)\| =0, t\ge t_k.
\end{equation}

And from (27) it follows tends to zero in probability of the sequence  $\bar{v}\left(\| x(t)\| \right)$ for $k\to \infty $ for all $ t \ge t_k, y\in {\it \mathbf{Y}}, h\in {\it \mathbf{H}}, x\in {\it \mathbb{R}}^{m} .$

So, from the properties of the Lyapunov function, we conclude that the non-negative supermartingale $v_{k} (\xi (t_{k} ),\eta _{k}, x(t_{k} ))$ for $k\to +\infty $ tends to zero in probability for all realizations of the process $\xi $ and sequence $\eta _{k} $.

Further, the nonnegative bounded supermartingale has a bound with probability 1 \cite{LM1}. Based on Theorem 1 (inequality (14) for the usual mathematical expectation), we obtain the asymptotic stochastically stability on the whole of the system (1)-(3) by the definition 3 (see (11)). Theorem 2 is proven.

\textit{End of proof of Theorem 2.}

{\textbf{Theorem 3.}} Suppose that the conditions of Theorem 2 are satisfied, and the Lyapunov functions $\{ v_{k} \} , \{ a_{k} \} , k\ge 0$ satisfy the inequalities
\begin{equation} \label{eq28}
c_{1} \| x\| ^{2} \le v_{k} (y, h, x)\le c_{2} \| x\| ^{2} ,
\end{equation}
\begin{equation} \label{eq29}
c_{3} \| x\| ^{2} \le a_{k} (y, h, x)\le c_{4} \| x\| ^{2}
\end{equation}
for some $c_{i} >0, i=\overline{1, 4},$ for all $k\in \mathbb{N}, y\in {\it \mathbf{Y}}, h\in {\it \mathbf{H}}, x\in {\it \mathbb{R}}^{m} .$

Then the system of random structure (1)-(3) is asymptotically stable in the mean square

\textit{Proof of Theorem 3.}
Using the inequality (19) for $k=0 $, based on (28) it is easy to obtain an inequality
$$
{\it \mathbf{E}}\left\{\| x(t_{N+1} )\| ^{2} \right\}\le \frac{1}{c_{1} } {\it \mathbf{E}}\{ v_{N+1} (\xi (t_{N+1} ),\eta _{N+1} , x(t_{N+1} ))\} \le 
$$
\begin{equation} \label{eq30}
\le \frac{1}{c_{1} } {\it \mathbf{E}}\{ v_{0 } (\xi (0),\eta _{0 } , x)\} \le \frac{c_{2} }{c_{1} } \| x\| ^{2}
\end{equation}
for all $N\ge 0, x\in {\it \mathbb{R}}^{m} $ and the initial distributions of the random vector $\{ \xi (0), \eta _{0} \} $.

Hence, by definition 4 (see (12)), it follows a stability in the mean square of the system of random structure (1)-(3).

Using the inequalities (20), (28) and (29), it is possible to obtain an inequality
$$
\sum _{k=0 }^{N}{\it \mathbf{E}}\left\{\| x(t_{N+1} )\| ^{2} \right\} \le \frac{1}{c_{3} } \sum _{k=0}^{N}{\it \mathbf{E}} \{ a_{k} (\xi (t_{k} ),\eta _{k} , x(t_{k} ))\} \le 
$$
$$
\le \frac{1}{c_{3} } {\it \mathbf{E}}\{ v_{0} (\xi (0),\eta _{0} , x)\} \le \frac{c_{2} }{c_{3} } \| x\| ^{2} .
$$

This inequality guarantees the convergence of the series whose members are ${\it \mathbf{E}}\left\{\| x(t_{k} )\| ^{2} \right\}$ for any initial data $x(0)=x$ and initial distributions of the random vector $\{ \xi (0), \eta _{0} \} $.

Therefore,
\[\mathop{\lim }\limits_{k\to \infty } \mathop{\sup }\limits_{y\in \mathbf{Y}, h\in \mathbf{H}, t\in [t_{k-1},t_k)} {\it \mathbf{E}}\left\{\| x(t)\| ^{2} \right\}=0,\]
which proves Theorem 3.

\textit{End of proof of Theorem 3.}

{\textbf{Theorem 4.}} [Corollary]
If the conditions of Theorem 3 are fulfilled and the inequality (28) holds, then the system of random structure (1)-(3) is stable in the mean square on the whole.

\section{Computation of the weak infinitesimal operator}

Based on the method \cite{LM12}, we will obtain an expression for calculating the explicit form of the weak infinitesimal operator (WIO) based on the system (1)-(3), which plays the role of the Lyapunov operator.

Let $U(t,y,h,x)$ be such a scalar integral function, that the sequence 
$$
\{v_k(y,h,x)\equiv U(t_k,y,h,x), k\ge 0\}
$$
is a Lyapunov function.

It is possible to prove \cite{LM1} that the pair $(\xi(t),x(t),t\ge 0,)$ is a Markov process and it is possible to introduce WIO
$$
(\mathcal{L}U)(t,y,h,x):=\lim\limits_{\Delta t\downarrow0}\frac{1}{\Delta t}[\mathbf{E}_{y,h,x}^{(t)}\{U(t+\Delta t,\xi(t+\Delta t),\eta(t+\Delta t),x(t+\Delta t))-
$$
\begin{equation}\label{eq31}
-U(t,y,h,x)\} ], 
\end{equation}
where $\mathbf{E}_{y,h,x}^{(t)}U=\mathbf{E}\{U\vert \xi(t)=y,\eta(t)=h,x(t)=x\}$, $\eta(t):=\eta_k$ and $t_k\le t <t_{k+1}, k\ge 0$. It is natural to assume that the function $U$, defined above, belongs to the domain of definition of the operator $\mathcal{L}$, if the limit (31) exists in the sense of uniform convergence in some neighborhood of the point $(y,x)$ uniformly by $h\in \mathbf{H}$. 

Let's introduce the operator $\mathcal{L}_0$ which is related to Markov switchings (2) at the moment $t_k, k\ge 0$:
\begin{equation}\label{eq32}
(\mathcal{L}_0U):=\mathbb{I}_{t\in K}\left[\int\limits_{\mathbf{H}}U(t,y,h,x)\mathbf{P}_k(h,dz)-U(t,y,h,x)\right],
\end{equation}
where $\mathbf{P}_k(h,dz)$ is the transition probability of the Markov chain at the $k$-th step, $\mathbb{I}$ is the indicator of the set $K$.

At the moment $\tau$ of changing of the structure of the parameter  $\xi$ of the system $y_i\rightarrow y_j$ there is a jump-like change in the phase vector $x$ with transition probability
\begin{equation}\label{eq33}
p_{ij}(\tau,x,A):=\mathbf{P}\{x(\tau)\in A \vert x(\tau-)=x, \xi(\tau-)=y_i, \xi(\tau)=y_j\}, A\subset \mathbb{R}^m.
\end{equation}

{\textbf{Theorem 5.}}
The weak infinitesimal operator $\mathcal{L}$ on the solutions of the system (1)-(3) of the function $U$ is calculated by the formula
\begin{equation}\label{eq34}
(\mathcal{L}U)(t,y,h,x)=(\mathcal{L}_tU)(t,y,h,x)+(\mathcal{L}_xU)(t,y,h,x)+(\mathcal{L}_yU)(t,y,h,x)+(\mathcal{L}_0U)(t,y,h,x),
\end{equation}
where
\begin{equation}\label{eq35}
(\mathcal{L}_tU)(t,y,h,x)=\frac{\partial U(t,y,h,x)}{\partial t},
\end{equation}
\begin{equation}\label{eq36}
(\mathcal{L}_xU)(t,y,h,x)=(\nabla_xU,a(t,y,x))+\frac{1}{2}Sp(\nabla^2_{xx}Ub(t,y,x),b^T(t,y,x)),
\end{equation}
\begin{equation}\label{eq37}
(\mathcal{L}_yU)(t,y,h,x)=\sum\limits_{i\neq j}\left[\int\limits_{\mathbb{R}^m}U(t,y_j,h,\zeta)p_{ij}(t,x,d\zeta)-U(t,y_i,h,x)\right]q_{ij}.
\end{equation}
Here $(\cdot,\cdot)$ is scalar product; $(\Delta U)=(\frac{\partial U}{\partial x_1},...,\frac{\partial U}{\partial x_m})^T$, $\frac{\partial U}{\partial x_i}, i={1,...,m}$ is the derivative of the $i$-th coordinate of the vector $x \in \mathbb{R}^m$; $\nabla^2_{xx}U=\left[\frac{\partial^2 U}{\partial x_i\partial x_j}\right]_{i,j=1}^m$ is a matrix of second derivatives; $Sp$ is a trace of the matrix; $q_{ij}=-\frac{\tilde{q}_{ij}}{\tilde{q}_i}$; $(\mathcal{L}_0U)(t,y,h,x)$ calculated by formula (32); $U$ is a function differentiable with respect to $t$, which has derivatives of the 1st and 2nd order by the last argument.

\textit{Proof of Theorem 5.}
By definition (31)
$$
(\mathcal{L}U)(t,y,h,x):=\lim\limits_{\Delta t\downarrow0}\frac{1}{\Delta t}[\mathbf{E}_{y,h,x}^{(t)}\{U(t+\Delta t,\xi(t+\Delta t),\eta(t+\Delta t),x(t+\Delta t))-
$$
$$
-U(t,y,h,x)\} ].
$$
Next,
$$
(\mathcal{L}U)(t,y,h,x):=\lim\limits_{\Delta t\downarrow0}\frac{1}{\Delta t}[\mathbf{E}_{y,h,x}^{(t)}\{U(t+\Delta t,\xi(t+\Delta t),\eta(t+\Delta t),x(t+\Delta t))-
$$
$$
-U(t,y,h,x)\pm U(t,\xi(t+\Delta t),\eta(t+\Delta t),x(t+\Delta t)) \pm
$$
$$
\pm U(t,y,\eta(t+\Delta t),x(t+\Delta t))\pm U(t,y,h,x(t+\Delta t)) \}].
$$

Therefore, $\mathcal{L}$ can be represented as
$$
(\mathcal{L}U)(t,y,h,x):=\lim\limits_{\Delta t\downarrow0}\frac{1}{\Delta t}[\mathbf{E}_{y,h,x}^{(t)}\{U(t+\Delta t,\xi(t+\Delta t),\eta(t+\Delta t),x(t+\Delta t))-
$$
$$
-U(t,\xi(t+\Delta t),\eta(t+\Delta t),x(t+\Delta t))\}]+
$$
$$
+\lim\limits_{\Delta t\downarrow0}\frac{1}{\Delta t}[\mathbf{E}_{y,h,x}^{(t)}\{U(t,\xi(t+\Delta t),\eta(t+\Delta t),x(t+\Delta t))-
$$
$$
-U(t,y,\eta(t+\Delta t),x(t+\Delta t))\}]+
$$
$$
+\lim\limits_{\Delta t\downarrow0}\frac{1}{\Delta t}[\mathbf{E}_{y,h,x}^{(t)}\{U(t,y,\eta(t+\Delta t),x(t+\Delta t))-
$$
$$
-U(t,y,h,x(t+\Delta t))\}]+
$$
$$
+\lim\limits_{\Delta t\downarrow0}\frac{1}{\Delta t}[\mathbf{E}_{y,h,x}^{(t)}\{U(t,y,h,x(t+\Delta t))-U(t,y,h,x\}].
$$

Let's consider each term separately. 

The form of the first term $\mathcal{L}_tU$ is obvious.

Let's establish the explicit form of the term $\mathcal{L}_xU$. Consider a complete group of disjoint events constructed as follows: denote by $H_i$ the event which means that the structure (1) does not change in the interval $(t,t+\Delta t]$, i.e. $\xi(\tau) =y_i$ at $\tau\in(t,t+\Delta t]$. Then with an accuracy of $o(\Delta t)$ we obtain \cite{LM13}
$$
\mathbf{P}(H_i)=-q_i\Delta t.
$$

Next, denote by $H_{ij}$ event, which means that in the interval $(t,t+\Delta t]$ a change $y_i\rightarrow y_j\neq y_i$ occurs. Then with accuracy up to $o(\Delta t)$ we have
$$
\mathbf{P}(H_{ij})=-q_{ij}\Delta t.
$$

Denote by $\Delta_iU:=U(t+\Delta t,\xi(t+\Delta t),h,x(t+\Delta t)-U(t,y_i,h,x)$ and by $\Delta_{ij}U$ the increment $\Delta U$ upon occurrence of the event $H_{ij}$. Let's calculate the increments $\Delta_iU$ and $\Delta_{ij}U$ of the function $U$ when events $H_i, H_{ij}, i\neq j$, occur, neglecting terms of order $o(\Delta t)$:
\begin{equation}\label{eq38}
\Delta_iH=\left[\frac{\partial U}{\partial t}+(\nabla_xU,a(t,y_i,x))+\frac{1}{2}Sp(\nabla^2_{xx}U b(t,y_i,x), b^T(t,y_i,x))\right]\Delta t+o(\Delta t).
\end{equation}
Here the partial derivatives are calculated at a point $(t,y_i,x)$, where $x$ is the solution of equation (1) with initial condition $\xi(t)=y_i, x(t)=x, s>t\ge0$. Next, for $\mathcal{L}_yU$ in the case of a change in the structure $y_i\rightarrow y_j$ in the interval $(t.t+\Delta t]$ we will get an increase 
\begin{equation}\label{eq39}
\Delta_{ij}U=U(t+\Delta t, y_j,h,x(t+\Delta t))-U(t,y_i,h,x)
\end{equation}
with the probability $q_{ij}\Delta t$.

The terms that illustrate the possibility of changing the structure of $\xi$ are not included in the last equality and there are no Markov switchings, since after averaging they have the order of $o(\Delta t)$ and we can ignore them.

To calculate $\mathbf{E}\{\Delta U\vert \xi(t)=y_i,\eta_0=h,x(t)=x\}$, we use the full probability formula
$$
\mathbf{E}\{\Delta U\vert \xi(t-)=y_i,\eta_0=h,x(t)=x\}=
$$
$$
=\mathbf{E}\{\mathbf{E}\{\Delta U\vert \xi(t)=y_j,\xi(t-)=y_i,\eta_0=h,x(t)=x\}\},
$$
where the external mathematical expectation on the right-hand side is calculated by the variable $\xi$ at the moment $t$.

Ignoring terms of order $o(\Delta t)$, from (38) and (39) we obtain
$$
\mathbf{E}\{\Delta U\vert \xi(t)=y_i,\eta_0=h,x(t)=x\}=
$$
$$
=\left[\frac{\partial U}{\partial t}+(\nabla_xU,a(t,y_i,x))+\frac{1}{2}Sp(\nabla^2_{xx}U b(t,y_i,x), b^T(t,y_i,x))\right](1-q_i\Delta t)\Delta t+
$$
$$
+\sum\limits_{i\neq j}\left[\int\limits_{\mathbb{R}^m}U(t,y_j,h,\zeta)p_{ij}(t,x,d\zeta)-U(t,y_i,h,x)\right]q_{ij}\Delta t+o(\Delta t).
$$

When calculating the third term, we used the property $x^TBx=Sp(Bxx^T)$ and the property of the Wiener process with respect to the covariance of the increment \cite{LM1}, \cite{LM9}.

Using division by $\Delta t$ and passing to the boundary at $\Delta \downarrow 0$, we obtain the first, second, and third terms in (34). The idea of calculating the fourth term $\mathcal{L}_0U$ can be found in \cite{LM7}, pp. 163-164. Theorem 5 is proved.

\textit{End of proof of Theorem 5.}

\section{Stability in probability on the whole of a linear stochastic system of random structure}

One-dimensional linear stochastic system of random structure given by SDE
\begin{equation}\label{eq40}
dx(t)=a(\xi(t))x(t)dt+b(\xi(t)),x(t)dw(t),~t\in \mathbb{R}_{+}\backslash K,
\end{equation}
with Markov switchings
$$
\Delta x(t)=g(t_k -,\xi(t_k -),\eta_k,x(t_k -)),
$$
\begin{equation}\label{eq41}
t_k \in K=\{t_n \Uparrow\}, \lim\limits_{n\rightarrow\infty}t_n=t^{*}\in [0,T<\infty],
\end{equation}
and initial conditions
\begin{equation}\label{eq42}
x(0)=x_0\in \mathbb{R}^1,~\xi(0)=y\in\mathbf{Y}, \eta_{0}=h\in\mathbf{H}, 
\end{equation}
where $x\in \mathbb{R}^1$ is a strong solution of SDE (40); $t^*$ is a concentration point; $\xi$ is a Markov chain with a finite number of states $\mathbf{Y}=\{1,2,...,N_{\xi} \}$ and generator $Q=\{\title{q}_{ij}\}, i,j=\{1,...,N_{\xi} \}$; $\{\eta_k, k\geq0\}$ is a Markov chain with values in space $\mathbf{H}$ and the transition probability at the $k$-th step $\mathbb{P}_k(h,dz)$; $w(t), t\ge0$ is a one-dimensional standard Wiener process; the processes $w, \xi$ and $\eta$ are independent \cite{LM1}, \cite{LM2}.

We obtain sufficient conditions for the stability of the system (40)-(42) in probability on the whole.

Let's choose a Lyapunov function in the form \cite{LM13}
\begin{equation}\label{eq43}
v(\xi(t),h,x)=\gamma \xi(t)\| x\|^\beta, \gamma>0.
\end{equation}

Let the functions $a(i)=a_i$, $b(i)=b_i$ be such that for all $i=\{1,...,N_{\xi} \}$
\begin{equation}\label{eq44}
a_i-\frac{b^2_i}{2}<-\varepsilon.
\end{equation}

Then in (43)
\begin{equation}\label{eq45}
\beta=\varepsilon b^{-2}, b=\max\limits_{i=\{1,...,N\}}\{b_i\}.
\end{equation}

We can show which restrictions must satisfy the transitional probabilities $q_{ij}$ of the Markov chain $\xi$ and $\mathbb{P}_k(h,dz)$ of the Markov chain $\eta$, so that the system (40)-(42)  is stable in probability on the whole.

Calculating $\mathcal{L}v$ on the solutions of the system (40)-(42), we obtain
$$
(\mathcal{L}v)(t_k,y,h,x)=\gamma \| x\|^\beta \left\{ b_i\left(a_i+\frac{\beta-1}{2}b_i \right)+\sum\limits_{j\neq i}^{k}(j-i)q_{ij}\right\}+
$$
$$
+\int\limits_{\mathbf{H}}\gamma i \| x+g(t_k,y,h,x) \|^\beta \mathbb{P}_k(h,dz)-\gamma i \| x\|^\beta.
$$

Considering (43)-(45), at the point $(\xi(t)=i,x)$ we have 
$$
(\mathcal{L}v)(t_k,y,h,x)=\gamma \| x\|^\beta \left[ -\frac{\beta i\varepsilon}{2}+a_i  \right] +
$$
\begin{equation}\label{eq46}
+i\gamma \left[ \int\limits_{\mathbf{H}} \| x+g(t_k,y,h,x) \|^\beta \mathbb{P}_k(h,dz)- \| x\|^\beta \right],
\end{equation}
where $a_i=\sum\limits_{j>i}^{k}(j-i)q_{ij}, a_k=0$.

Assuming that for $\forall h\in \mathbf{H}$ of the Markov chain $\eta$ the transition probability at the $k$-th step $\mathbb{P}_k(h,dz)$ such that
\begin{equation}\label{eq47}
 \int\limits_{\mathbf{H}} \| x+g(t_k,y,h,x) \|^\beta \mathbb{P}_k(h,dz) \le 2\| x \|^\beta,
\end{equation}
then the right-hand side of (46) will take the form
$$
(\mathcal{L}v)(t_k,y,h,x)=\gamma \| x\|^\beta \left[ -\frac{\beta i\varepsilon}{2}+a_i +i \right]=\gamma \| x\|^\beta \left[ -\frac{ i(\beta\varepsilon+2)}{2}+a_i  \right].
$$

The function (43) satisfies the condition $\mathcal{L}v<0$ if the expression in square brackets is negative. Thus, we can formulate the following statement.

{\textbf{Theorem 6.}}
If the conditions  (44), (45) are met and 
\begin{equation}\label{eq48}
a_i<\frac{i(\beta\varepsilon+2)}{2}, i=\{1,...,N_\xi \},
\end{equation}
then the solution of the system (40)-(42) is stable in probability on the whole for all fixed $y\in\mathbf{Y}$ and $h\in\mathbf{H}$.

\section{Model example}

Consider the linear stochastic differential equation
\begin{equation}\label{eq49}
dx(t)=a(\xi(t))x(t)dt+b(\xi(t))x(t)dw(t), t\geq0, r>0,
\end{equation}
with impulse action
\begin{equation}\label{eq50}
\Delta x\left(2-\frac{1}{k}\right)=x\left(2-\frac{1}{k}-\right)+e^{-\alpha k\eta_k}\left(x\left(2-\frac{1}{k}\right)\wedge 1\right), k\rightarrow \infty,
\end{equation}
and initial condition
\begin{equation}\label{eq51}
x(0)=10, \xi(0)=y_0\in\mathbf{Y}, \eta_0=1.
\end{equation}

Here $a$ and $b$ are constants that depend on Markov process $\xi$ with values in dimensional space  $(\mathbf{Y},\mathcal{Y})$ with generator $Q$, and $\eta_k, k\geq0,$ is Markov chain with two non-absorbing states $h_1=0$ and $h_2=1$.

According to \cite{LM3} the solution of the system (49)-(51) exists, for example, when $\alpha=1.673$. 

Case 1. Let's consider the same coefficients as in \cite{LM3}:

- if $\xi=1$: $a=1, b=0.3$;

- if $\xi=2$: $a=-0.5, b=2.1$;

- $\eta_k\in\{1,2\}$.

In this case condition (44) is not hold for $i=1$ because 
$$
1-\frac{0.3^2}{2}=0.955>0.
$$

Therefore, the solution can be unstable. Indeed, if we consider an example of the realization of the solution of the system (49)-(51) with indicated parameters, then we observe a rapid growth (see Figure 1 (a)).

\begin{figure}[ht]
\centering
\includegraphics[width=0.9\textwidth]{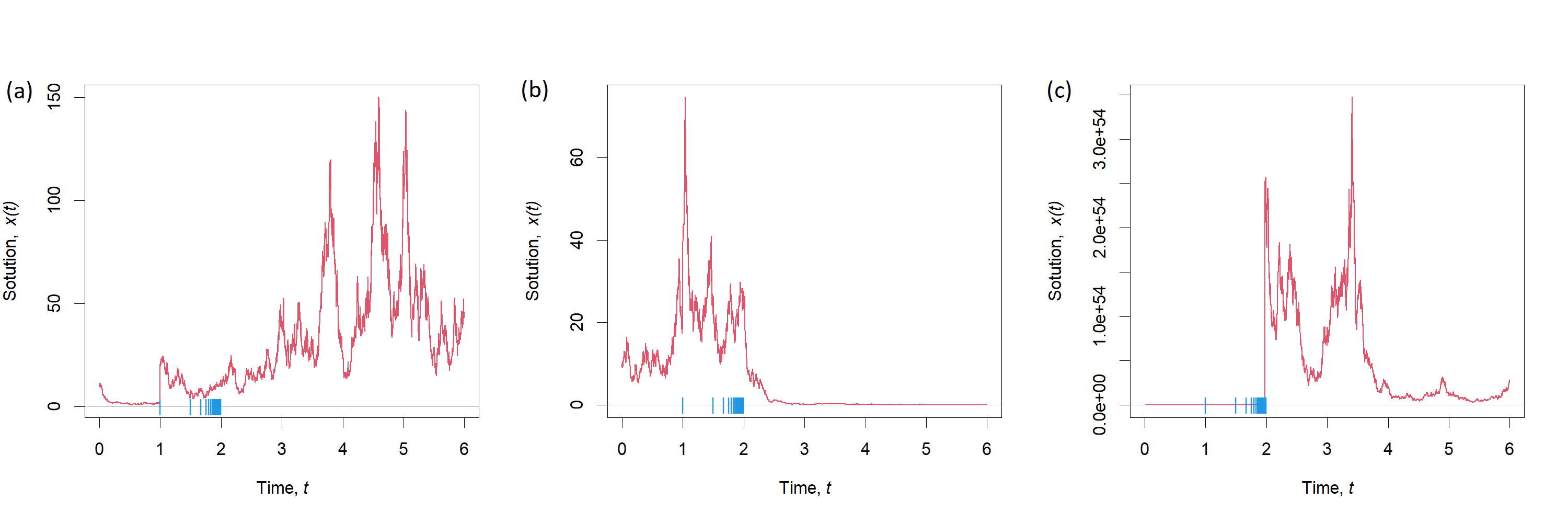}
\caption{ Solution trajectories: (a) case 1 (unstable), (b) case 2 (stable), (c) case 3 (unstable with an extreme growth at $t=2$). The red line corresponds to the system's solution $x(t)$ evolution. Blue marks - moments of impulse actions. 
}
\label{fig1}
\end{figure}

Case 2. Next, we consider the values of the coefficients:

- if $\xi=1$: $a=-1, b=0.3$;

- if $\xi=2$: $a=0.5, b=2$;

- $\eta_k\in\{1,2\}$.

Condition (44) for $i=1$ has the next form

$$
-1-\frac{0.3^2}{2}=-1.045<-\varepsilon,
$$
and for  $i=2$ has the form
$$
0.5-\frac{2^2}{2}=-1.5<-\varepsilon,
$$
and holds for $\varepsilon=0.1$.

According to (45)

$$
\beta=\frac{0.1}{2^2}=0.025.
$$

And (48) hold:

- if $i=1$: $-1<\frac{1\cdot(0.025\cdot0.1+2)}{2}=1.00125$;

- if $\xi=2$: $0.5<\frac{2\cdot(0.025\cdot0.1+2)}{2}=2.0025$.

So, all conditions of Theorem 6 are held and the solution of the system (49)-(51) with indicated parameters is stable in probability on the whole. Indeed, in the realization (see Figure 1 (b)) we observe a direction to zero after the point $t=3$.

Case 3. Here the values of the coefficients are the same as in case 2, but the impulse action has the next form
$$
\Delta x\left(2-\frac{1}{k}\right)=x\left(2-\frac{1}{k}-\right)+e^{\alpha k\eta_k}\left(x\left(2-\frac{1}{k}\right)\wedge 1\right), k\rightarrow \infty.
$$

In this case condition (47) does not hold and we cannot guarantee stability in the probability of solution of the system (49)-(51): we observe a very rapid growth (see Figure 1 (c)).

\section{Discussion}\label{sec12}

The condition (44) means that stability in probability can be ensured due to larger values of the coefficients and the fulfillment of the condition (47) even when the system is unstable
$$
dx(t)=a_ix(t)dt.
$$

For example, if we consider the second case of the model example with the coefficients $a=0.5, b=2$, we will see that the solution of the system corresponding to the deterministic part is not Lyapunov stable, but the solution of the stochastic system, as was demonstrated, is stable in probability.

\section{Conclusion}\label{sec13}

In this paper, we explore sufficient conditions for the asymptotic stability of stochastic differential equations with a random structure, particularly in the context of jump condensation points. Our main result is presented in Theorem 2, which leverages the second Lyapunov method and involves the construction of corresponding Lyapunov functions. An important consideration in analyzing systems of random structure is the relationship between the magnitudes of jumps, denoted as $L_k$, and the jump times, denoted as $\tau_k$. The implications of Theorem 2 are demonstrated through an example system whose stability can be modulated by varying parameters. We also highlight a remarkable observation that the system can maintain asymptotic stability even if, for some fixed value of the random process $\xi(t)$, the system described by equation (1) becomes unstable when jumps (2) are absent.

In future studies, we plan to investigate the stability of stochastic differential equations with a random structure, particularly when the jump moments, denoted as $\tau_k$, are random variables satisfying the condition
$$
P\left(\lim\limits_{k\rightarrow\infty}t_k=t_\infty<\infty\right)>0,
$$
This implies a non-zero occurrence of condensation points. Furthermore, the weak independence between the jumps and their corresponding moments will also be considered as part of this analysis.

\section*{Acknowledgments}
We would like to acknowledge the administrations of the Luxembourg Institute of Health (LIH) and Luxembourg National Research Fund (FNR) for their support in organizing scientific contacts between research groups in Luxembourg and Ukraine.

\end{document}